\documentclass{amsart}

\usepackage{amsfonts, amsmath, amsthm, amssymb, epsfig,
bm, graphics, psfrag, latexsym, chngcntr}
\usepackage[all,cmtip]{xy}

\newtheorem{thm}{Theorem}[section]         \newtheorem*{thm*}{Theorem}
                 \newtheorem*{lem*}{Lemma}
   \newtheorem*{prop*}{Proposition}
         \newtheorem*{ques*}{Question}
 \newtheorem*{defn*}{Definition}
   \newtheorem*{rem*}{Remark}

\newcommand{\wh}{\widehat}
\newcommand{\del}{\partial}
\newcommand{\Del}{\Delta}
\newcommand{\wt}{\widetilde}
\newcommand{\mc}{\mathcal}
\newcommand{\mb}{\mathbb}
\newcommand{\mf}{\mathfrak}

\newcommand{\Zoltan}{Zolt\'{a}n}
\newcommand{\Szabo}{Szab\'{o}}
\newcommand{\Ozsvath}{Ozsv\'{a}th}
\newcommand{\Andras}{Andr\'{a}s}
\newcommand{\Juhasz}{Juh\'{a}sz}
\newcommand{\Poincare}{Poincar\'{e}}

\begin{document}

\title{Maslov index formulas for Whitney $n$-gons}
\author{Sucharit Sarkar}
\address{Department of Mathematics, Columbia University, New York, NY 10027, USA}
\email{sucharit@math.columbia.edu}


\date{}

\begin{abstract}
In this short article, we find an explicit formula for Maslov index of Whitney $n$-gons joining intersections points of $n$ half-dimensional tori in the symmetric product of a surface. The method also yields a formula for the intersection number of such an $n$-gon with the fat diagonal in the symmetric product. 
\end{abstract}

\maketitle

\section{Introduction}\label{sec:introduction}

In \cite{POZSz3manifolds}, Peter \Ozsvath{} and \Zoltan{} \Szabo{} introduced Heegaard Floer homology, a powerful collection of invariants for closed oriented $3$-manifolds. There are various versions, but all of them involve counting the number of points in the unparametrized moduli space, coming from Maslov index one disks joining some intersection points of two half-dimensional tori in a certain symmetric product. It is then clear that, to achieve a combinatorial understanding of the theory, we need a formula for the Maslov index of such disks. In \cite{JR}, Jacob Rasmussen gave a formula, depending only on the combinatorial information coming from a $2$-chain in the Heegaard diagram representing such a disk, which relates the intersection number of the disk with the fat diagonal in the symmetric product, with its Maslov index. Robert Lipshitz gave a cylindrical reformulation of the whole theory in \cite{RL}, and using his reformulation, he was able to determine the Maslov index of such disks.

However in \cite{POZSz4manifolds}, \Ozsvath{} and \Szabo{} showed that given a cobordism between two $3$-manifolds, there is an induced map on the Heegaard Floer homologies. This converts Heegaard Floer homology into a $3+1$ dimensional TQFT, and leads to the definition of a smooth $4$-manifold invariant called the \Ozsvath-\Szabo{} invariant, which is conjecturally equal to the Seiberg-Witten invariant. This theory also involves counting the number of points in the moduli space, coming from Maslov index zero triangles joining three intersection points in three half-dimensional tori in some symmetric product. We will derive a formula to determine the Maslov index of such triangles. Our method will also yield a formula for the Maslov index of $n$-gons joining intersection points of $n$ half-dimensional tori in the symmetric product. Using the same ideas, we will also be able to compute the intersection number of such $n$-gons with the fat diagonal.

The outline of this rather short paper is as follows. After this introductory section, we proceed onto Section \ref{sec:setting}, where we will give the general setting for the whole construction, fix many notations, and explicitly write down our formula. In Section \ref{sec:properties}, we will state and prove many properties of the Maslov index and of our formula. In Section \ref{sec:main}, we will prove the main result that our formula actually does compute the Maslov index; and finally in Section \ref{sec:applications}, we will give some applications.

\subsection*{Acknowledgements}
I would like to thank Professor \Zoltan{} \Szabo{} for
suggesting this problem to me, for providing a general structure of the answer, for guiding me through the solution and for pointing out the application. I would also like to thank \Andras{} \Juhasz{} for having many helpful discussions regarding various definitions of the Maslov index. Finally, I would like to thank the referee for the careful review and the helpful comments. A part of the work was done when the author was partially supported by the Princeton centennial fellowship and when the paper finally took form, he was fully supported by the Clay postdoctoral fellowship. 

\section{Setting}\label{sec:setting}

The definition of a Heegaard diagram in its full generality can be bit overwhelming at first, but we ask the patient reader to hold tight for a little while. A Heegaard diagram is essentially a closed oriented surface $\Sigma$ with a collection of simple closed curves $\eta^i_j$. The various intersections among the $\eta$ curves are always assumed to be transverse. Throughout this article, the genus of the surface will be $g$; the parameter $i$ will range from $1$ to $n$, $n$ being the number of half-dimensional tori and the parameter $j$ will range from $1$ to $k$, $k$ being the dimension of each half-dimensional torus. Usually we have $k\geq g$, and in that case there are $(k-g+1)$ marked points on the surface $w_1,\ldots,w_{k-g+1}$. For each $i$, let $\eta^i=\{\eta^i_1,\ldots,\eta^i_k\}$ and let $w=\{w_1,\ldots,w_{k-g+1}\}$. Furthermore we assume that for each $i$, $\eta^i$ is a disjoint collection of curves, and if $k\geq g$, we assume that $\Sigma\setminus(\cup_j\eta^i_j)$ has $(k-g+1)$ components (i.e. $\eta^i$ spans a half-dimensional subspace of $H_1(\Sigma)$) each containing some basepoint $w_l$. A Heegaard diagram $\mc{H}=(\Sigma,\eta^1,\ldots,\eta^n,w)$ encapsulates this whole structure.

Heegaard Floer theory usually just deals with small values of $n$, and up to $n=3$, we use the Greek letters $\alpha$, $\beta$ and $\gamma$ to denote $\eta^1$, $\eta^2$ and $\eta^3$ respectively. For $n=1$ and $k\geq g$, the Heegaard diagram $(\Sigma,\alpha,w)$ denotes a genus $g$ handlebody $U_{\alpha}$ obtained by first thickening $\Sigma$ to $\Sigma\times[0,1]$, then adding $k$ two-handles along $\alpha_j\times\{1\}$, and finally adding $(k-g+1)$ three-handles to the $(k-g+1)$ boundary components each marked with a point $w_l\times\{1\}$. For $n=2$ and $k\geq g$, the Heegaard diagram $(\Sigma,\alpha,\beta,w)$ represents a closed oriented $3$-manifold $Y_{\alpha,\beta}$ obtained by gluing together $U_{\alpha}$ and $-U_{\beta}$. For $n=3$ and $k\geq g$, the Heegaard diagram represents a smooth $4$-manifold $W_{\alpha,\beta,\gamma}$ with $\del W_{\alpha,\beta,\gamma}=Y_{\alpha,\beta}\cup Y_{\beta,\gamma}\cup Y_{\gamma,\alpha}$, obtained by first taking a triangle $\Delta ABC$, and then gluing together $\Delta ABC\times\Sigma$, $AB\times U_{\gamma}$, $BC\times U_{\alpha}$ and $CA\times U_{\beta}$.

Given a Heegaard diagram, all the Heegaard Floer invariants are constructed in essentially the same way. The ambient manifold is the symmetric product $Sym^k\Sigma=(\prod^k_{j=1}\Sigma)/S_k$ and let $T_i=\prod_j \eta^i_j$ be $n$ totally real half-dimensional tori lying inside the symmetric product. Let $P^{i,j}=T_i\cap T_j$ and if $p\in P^{i,j}$, then $p$ is an unordered $k$-tuple of points, each lying on $\Sigma$. We call those points the coordinates of $p$ and write $p=(p_{1},\ldots,p_k)$. The basic algebraic object that we study is the module over some ring (usually $\mb{Z}$ or $\mb{Z}[U_1,\ldots,U_{k-g+1}]$) freely generated by the finitely many points in $\cup_{i,j}P^{i,j}$.

This is the point from where the story gets complicated. If $n\geq 2$, for each value of $i$, we choose a point $p^{i,i+i}$ from $P^{i,i+1}$ (the counting throughout being done modulo $n$). Let $D^2$ be the unit disk in the complex plane with $n$ fixed marked points on the boundary numbered in a counter-clockwise fashion as $t_1,t_2,\ldots,t_n$ and let $s_i\subset \del D^2$ be the positively oriented arc joining $t_{i-1}$ to $t_i$. We consider maps from $D^2$ to the symmetric product which map $t_i$ to $p^{i,i+1}$ and $s_i$ to $T_i$ for all values of $i$. A Whitney $n$-gon is a homotopy type of such maps, where the homotopy is through maps respecting the same boundary conditions. Let $\pi_2(p^{1,2},\ldots,p^{n,1})$ denote the set of all Whitney $n$-gons connecting the points $p^{1,2},\ldots,p^{n,1}$. If $\phi\in\pi_2(p^{1,2},\ldots,p^{n,1})$ is a Whitney $n$-gon, let $\iota(\phi)$ be the intersection number of $\phi$ with the fat diagonal in the symmetric product. Note that this is well-defined since the boundary of $\phi$ lies in $\cup_i T_i$ and the fat diagonal is disjoint from $\cup_i T_i$.

\begin{figure}[ht]
\psfrag{t}{$t^{*}$}
\psfrag{t12}{$t_1$}
\psfrag{s12}{$t_1$}
\psfrag{t23}{$t_2$}
\psfrag{s23}{$t_2$}
\psfrag{t34}{$t_3$}
\psfrag{s34}{$t_3$}
\psfrag{t45}{$t_4$}
\psfrag{s45}{$t_4$}
\psfrag{t51}{$t_5$}
\psfrag{s51}{$t_5$}
\begin{center}
\includegraphics[width=250pt]{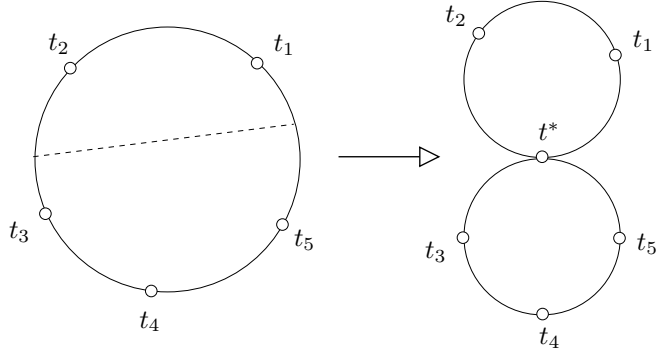}
\end{center}
\caption{The map that induces the addition operation on the Whitney $n$-gons}\label{fig:addition}
\end{figure}

The set of all Whitney $n$-gons enjoys a very nice multiplicative structure, which is worth mentioning at this point. As in the previous paragraph, let us choose points $p^{i,i+1}\in P^{i,i+1}$ for all values of $i$, and let us also choose an additional point $q\in P^{1,m}$ for some $1\leq m\leq n$. If we ignore the tori $T_i$ for $m<i\leq n$, we can talk about the Whitney $m$-gons $\pi_2(p^{1,2},\ldots,p^{m-1,m},q)$; similarly if we ignore the tori $T_i$ for $1<i<m$, we can talk about the Whitney $(n-m+2)$-gons $\pi_2(q,p^{m,m+1},\ldots,p^{n-1,n},p^{n,1})$. There is a map from $\pi_2(p^{1,2},\ldots,p^{m-1,m},q)\times \pi_2(q,p^{m,m+1},\ldots,p^{n-1,n},p^{n,1})$ to $\pi_2(p^{1,2},\ldots,p^{n,1})$, induced from the map from a disk with $n$ marked points on its boundary to a wedge of two disks with $m$ and $(n-m+2)$ marked points on their boundaries respectively. Figure \ref{fig:addition} illustrates such a map for $n=5$ and $m=3$ obtained by quotienting out the dotted arc to a point. The points marked $t_i$ map to $p^{i,i+1}$ in the symmetric product; the point marked $t^{*}$ along which the wedge is taken, maps to the point $q$. This map is denoted by the symbol $*$, i.e. given a Whitney $m$-gon $\phi\in \pi_2(p^{1,2},\ldots,p^{m-1,m},q)$ and a Whitney $(n-m+2)$-gon $\psi\in \pi_2(q,p^{m,m+1},\ldots,p^{n-1,n},p^{n,1})$, we get a Whitney $n$-gon $(\phi *\psi)\in\pi_2(p^{1,2},\ldots,p^{n,1})$.

The combinatorial nature of the theory ends here. A complex structure is chosen on the symmetric product $Sym^k(\Sigma)$, which is a generic perturbation of the complex structure induced from some complex structure on $\Sigma$. If $\phi\in\pi_2(p^{1,2},\ldots,p^{n,1})$ is a Whitney $n$-gon, then let $\mc{M}(\phi)$ be the moduli space of holomorphic maps from $D^2$ to $Sym^k(\Sigma)$ that represent $\phi$, and let the Maslov index $\mu(\phi)$ be the expected dimension of $\mc{M}(\phi)$. Note that here we fix the complex structure on $D^2\setminus\{t_1,\ldots,t_n\}$. If instead, for $n>3$ we allow the complex structure of the source to vary and for $n<3$ we quotient out the moduli space by the action of precomposition by automorphisms of the source, we get a very different moduli space $\wt{\mc{M}}(\phi)$ whose expected dimension is $\mu(\phi)+(n-3)$. All the Heegaard Floer invariants require the count of the number of points in $\wt{\mc{M}}(\phi)$ when its expected dimension is zero, i.e. when $\mu(\phi)=3-n$. Therefore given a Whitney $n$-gon $\phi$, it is an important problem to determine its Maslov index $\mu(\phi)$.

There is an alternate definition of the Maslov index $\mu$ that merits some attention. Let $\mc{R}_k$ be the space of all $k$-dimensional totally real subspaces of $\mb{C}^k$. There is a canonical generator $[G]\in H_1(\mc{R}_k)=\mb{Z}$. For every $i,j$ and for every point $p\in P^{i,j}$, choose a path $\tau(p)$ in $\mc{R}_k$ joining $T_p(T_i)\subset T_p(Sym^k(\Sigma))=\mb{C}^k$ to $T_p(T_j)\subset T_p(Sym^k(\Sigma))=\mb{C}^k$ which has index zero \cite{RS}. For $n\geq 2$, given a Whitney $n$-gon $\phi\in\pi_2(p^{1,2},\ldots,p^{n,1})$, we take the pullback of the tangent bundle of $Sym^k(\Sigma)$ to an $k$-dimensional complex bundle over the unit disk $D^2$. The disk being contractible, it is the trivial bundle. However the arc $s_i\subset\del D^2$ maps to the tori $T_i$; the pullback of the tangent bundle of $T_i$ produces a path $\tau'(s_i)$ in $\mc{R}_n$.  Then the loop $\tau'(s_1)\tau(\phi(t_1))\tau'(s_2)\cdots\tau(\phi(t_n))$ represents the homology class $\mu(\phi)[G]$.

Whitney $n$-gons are often represented by their shadows on the Heegaard surface $\Sigma$. Before we describe what we mean by this, let us quickly set up a few more notations. A region is a component of $\Sigma\setminus(\cup_{i,j}\eta^i_j)$. If a region is topologically a disk, and has $m$ $\eta$-arcs on its boundary, then we call such a region an $m$-sided region. The usual convention is to call the $2$-sided regions bigons and the $3$-sided regions triangles; however, the same convention dictates that we call the Whitney $2$-gons bigons and the Whitney $3$-gons triangles. To avoid any confusion, henceforth we will not use the terms bigons or triangles to denote either of these objects. The objects in the symmetric product will always be referred to as Whitney $n$-gons, and the regions on the Heegaard surface will always be referred to as $m$-sided regions. 

Let $D$ be a $2$-chain generated by the regions over $\mb{Z}$, i.e. $D=\sum_R a(R) R$ for some integers $a(R)$. Let $\del_i(D)=\del D_{|\eta^i}$. Therefore $\del_i(D)$ is an $1$-chain lying on the $\eta^i$ curves whose endpoints are the various intersection points between the $\eta^i$ curves and the other $\eta$ curves. For $n\geq 2$, given points $p^{i,i+1}\in P^{i,i+1}$, a domain joining them is a $2$-chain $D$ such that $\del (\del_i D)=\sum_l (p^{i,i+1}_l-p^{i-1,i}_l)$ for all values of $i$. If $p\in P^{i,j}$, we often misuse notation to write $p=\sum_l p_l$  and thereby write $\del(\del_i D)=p^{i,i+1}-p^{i-1,i}$. The set of all such domains is denoted by $\mc{D}(p^{1,2},\ldots,p^{n,1})$. It is an useful fact to remember that if $n\geq 3$ and if $D\in\mc{D}(p^{1,2},\ldots,p^{n,1})$, then $\del(\del_i(D))_{|\eta^{i+1}}=p^{i,i+1}=\sum_l p^{i,i+1}_l$.

Given a region $R$, let us choose a point $r\in R$ and let us consider the divisor $Z_r=\{r\}\times Sym^{k-1}(\Sigma)$. The divisor $Z_r$ is disjoint from $\cup_i T_i$, therefore given a Whitney $n$-gon $\phi$, the intersection number $\phi\cdot Z_r$ is well-defined. After the immediate observation that the intersection number is independent of the choice of the point $r$, we denote it by $n_R(\phi)$ and call it the coefficient of $\phi$ at $R$. Let $D(\phi)= \sum_R n_R(\phi) R$ be the shadow of $\phi$. It is easy to see that $D(\phi)\in\mc{D}(p^{1,2},\ldots,p^{n,1})$ and it is called the domain representing $\phi$. The map from $\pi_2(p^{1,2},\ldots,p^{n,1})$ to $\mc{D}(p^{1,2},\ldots,p^{n,1})$ given by $\phi$ mapping to $D(\phi)$, happens to be a bijection for high enough values of $k$ and $g$. If $q\in P^{1,m}$, then the addition operation on $2$-chains gives a map from $\mc{D}(p^{1,2},\ldots,p^{m-1,m},q)\times \mc{D}(q,p^{m,m+1},\ldots,p^{n-1,n},p^{n,1})$ to $\mc{D}(p^{1,2},\ldots,p^{n,1})$, which corresponds to the multiplication map for Whitney $n$-gons.

Given an $1$-chain $a_i$ supported on $\eta^i$ subject to the condition that all of its boundary points lie on the intersections of various $\eta$ curves, and another $1$-chain $a_j$ supported on $\eta^j$ subject to the same condition, we can define $a_j\cdot a_i$ as follows. First, recall that we have fixed a complex structure, and thereby an orientation, on $\Sigma$. Now orient all the $\eta^j$ circles. This gives four possible directions to translate $a_j$, such that if $\wt{a_j}$ is a small translate in one of the four directions, then no end point of $a_i$ lies on $\wt{a_j}$, and no endpoint of $\wt{a_j}$ lies on $a_i$. Therefore using the orientation on $\Sigma$, we can define the intersection of each of the four translates with $a_i$; the intersection number $a_j\cdot a_i$ is defined to be the average of the four. This is clearly seen to be well-defined and skew-symmetric, i.e. $a_j\cdot a_i=-a_i\cdot a_j$. 

Let $p=(p_1,\ldots,p_n)\in\cup_{i,j} P^{i,j}$ and let $D$ be a $2$-chain on $\Sigma$. The coefficient of $D$ at $p_l$ is denoted by $\mu_{p_l}(D)$ and is defined to be the average of the coefficients of $D$ in the four regions around $p_l$. The point measure $\mu_p(D)$ is defined as $\mu_p(D)=\sum_l \mu_{p_l}(D)$.

For a $2$-chain $D$, the Euler measure $e(D)$ is defined as follows. We fix a metric on $\Sigma$ under which all the $\eta$ curves are geodesics and all the intersections among the $\eta$ curves are at right angles. Then $e(\phi)$ is $\frac{1}{2\pi}$ times the integral of the curvature along the $2$-chain $D$. Being an integral, the Euler measure is additive. For an $m$-sided region, the Euler measure is $(1-\frac{m}{4})$.

Amid this rather long and dry section on notations, we might have lost track of our original goal. For $n\geq 2$, given a Whitney $n$-gon $\phi\in\pi_2(p^{1,2},\ldots,p^{n,1})$, we are trying to find formulas for $\iota(\phi)$, the intersection number with the fat diagonal and $\mu(\phi)$, the Maslov index. Without further ado, we start with a domain $D\in\mc{D}(p^{1,2},\ldots,p^{n,1})$, we implicitly assume $n\geq 2$, we present our candidates
$$\iota(D)=\mu_{p^{n,1}}(D)+\mu_{p^{1,2}}(D)+\sum_{n\ge j>l>1}\del_j(D)\cdot\del_l(D)-e(D),$$
$$\mu(D)=\iota(D)+2e(D)-\frac{k(n-2)}{2},$$
and we proceed onto the next section.

\section{A few properties of $\mu(\phi)$, $\mu(D)$, $\iota(\phi)$ and $\iota(D)$}\label{sec:properties}

First, we prove an important theorem that relates the point measures of different points.

\begin{thm}\label{thm:mainproperty}
Let $D\in\mc{D}(p^{1,2},\ldots,p^{n,1})$ be a domain and let $D'$ be a (possibly different) $2$-chain in the Heegaard diagram $\mc{H}$. Then $\mu_{p^{i,i+1}}(D')-\mu_{p^{i-1,i}}(D')=\del D'\cdot\del_i(D)=\sum_{j\ne i} \del_j(D')\cdot\del_i(D)$.
\end{thm}

\begin{proof}
After fixing the orientations on $\Sigma$ and the $\eta^i$ circles, we can assume that the $\eta^i$ curves run in the north-south direction, and every other $\eta$ curve intersects them perpendicularly in an east-west direction. Therefore for each point on the $\eta^i$ circles, we have four well defined corners, north-east, north-west, south-west and south-east. We can translate $\del_i(D)$ slightly in each of these four directions to get $\del_i(D)_{NE}$, $\del_i(D)_{NW}$, $\del_i(D)_{SW}$ and $\del_i(D)_{SE}$ respectively.

By travelling along $\del_i(D)_{NE}$ from the north-east corner of $p^{i-1,i}$ to the north-east corner of $p^{i,i+1}$, we see that, (coefficient of $D'$ at the north-east corner of $p^{i,i+1}$)$=$(the coefficient of $D'$ at the north-east corner of $p^{i-1,i}$)$+\del D'\cdot\del_i(D)_{NE}$. We have similar results for the other corners, and after taking averages, we get
$\mu_{p^{i,i+1}}(D')-\mu_{p^{i-1,i}}(D')=\del D'\cdot\del_i(D)=\sum_j \del_j(D')\cdot\del_i(D)
=\sum_{j\ne i} \del_j(D')\cdot\del_i(D)$.
\end{proof}

Next, we show that formulas for $\mu(D)$ and $\iota(D)$ are cyclically symmetric in $p^{1,2},\ldots,p^{n,1}$. This is immediate from the following theorem.

\begin{thm}
If $D\in\mc{D}(p^{1,2},\ldots,p^{n,1})$, then the expression $\mu_{p^{n,1}}(D)+\mu_{p^{1,2}}(D)+\sum_{n\ge j>l>1}\del_j(D)\cdot\del_l(D)$ is a cyclically symmetric expression in $p^{1,2},\cdots,p^{n,1}$.
\end{thm}

\begin{proof}
We only need to show that $\sum_{n\ge j>l>1}\del_j(D)\cdot\del_l(D)+\mu_{p^{n,1}}(D)+\mu_{p^{1,2}}(D)
=\sum_{n>j>l\ge 1}\del_j(D)\cdot\del_l(D)+\mu_{p^{n-1,n}}(D)+\mu_{p^{n,1}}(D)$.
Theorem \ref{thm:mainproperty} implies that $\mu_{p^{1,2}}(D)=\mu_{p^{n,1}}(D)+\sum_{j\ne 1}\del_j(D)\cdot\del_1(D)$, therefore the expression $\mu_{p^{n,1}}(D)+\mu_{p^{1,2}}(D)+\sum_{n\ge j>l>1}\del_j(D)\cdot\del_l(D)$ can be rewritten as $2\mu_{p^{n,1}}(D)+\sum_{n\ge j>l\ge 1}\del_j(D)\cdot\del_l(D)$. Similarly, we have $\mu_{p^{n-1,n}}(D)=\mu_{p^{n,1}}(D)+\sum_{j\ne n}\del_n(D)\cdot\del_j(D)$, therefore the expression $\mu_{p^{n-1,n}}(D)+\mu_{p^{n,1}}(D)+\sum_{n>j>l\ge 1}\del_j(D)\cdot\del_l(D)$ can also be rewritten as $2\mu_{p^{n,1}}(D)+\sum_{n\ge j>l\ge 1}\del_j(D)\cdot\del_l(D)$. This finishes the proof.
\end{proof}

Staying with the same notations, let $p^{i,i+1}\in P^{i,i+1}$ and for some $1\le m\le n$, let $q\in P^{1,m}$. If $\phi\in\pi_2(p^{1,2},\ldots,p^{m-1,m},q)$ and $\psi\in \pi_2(q,p^{m,m+1},\ldots,p^{n-1,n},p^{n,1})$, the multiplication operation gives rise to $(\phi *\psi)\in\pi_2(p^{1,2},\ldots,p^{n,1})$. From the definition of the intersection number $\iota$ and from the alternate definition of the Maslov index $\mu$, it is easy to prove that $\iota(\phi *\psi)=\iota(\phi)+\iota(\psi)$ and $\mu(\phi *\psi)=\mu(\phi)+\mu(\psi)$. A very similar result holds for domains in the Heegaard diagram.

\begin{thm}\label{thm:addition}
  For $p^{i,i+1}\in P^{i,i+1}$ and $q\in P^{1,m}$, let   $D_1\in\mc{D}(p^{1,2},p^{2,3},\ldots,\allowbreak
p^{m-1,m},q)$ and   $D_2\in\mc{D}(q,p^{m,m+1},\ldots, p^{n-1,n},p^{n,1})$. Then   $\iota(D_1+D_2)=\iota(D_1)+\iota(D_2)$ and   $\mu(D_1+D_2)=\mu(D_1)+\mu(D_2)$.
\end{thm}

\begin{proof}
We know that $e(D_1+D_2)=e(D_1)+e(D_2)$ and $\frac{k(m-2)}{2}+\frac{k(n-m+2-2)}{2}=\frac{k(n-2)}{2}$, so we only need to prove
\begin{align*}
\mu_{p^{n,1}}(D_1+D_2)+\mu_{p^{1,2}}(D_1+D_2)+\sum_{n\ge j>l>1}\del_j(D_1+D_2)\cdot\del_l(D_1+D_2)\\
=\sum_{m\ge j>l>1}\del_j(D_1)\cdot\del_l(D_1)
+\sum_{n\ge j>l\ge m}\del_j(D_2)\cdot\del_l(D_2)\\
+\mu_{p^{1,2}}(D_1)+\mu_{p^{n,1}}(D_2)+\mu_q(D_1+D_2).
\end{align*}

Two straightforward applications of Theorem \ref{thm:mainproperty} imply that 
$$\mu_q(D_1)=\mu_{p^{n,1}}(D_1)-\del_1(D_2)\cdot\del D_1 {\text{ and}}$$
$$\mu_q(D_2)=\mu_{p^{1,2}}(D_2)-\del D_2\cdot\del_1(D_1).$$ Observe that $\del_j D_1=0$ for $n\ge j>m$ and $\del_j D_2=0$ for $m>j>1$. Therefore,
\begin{align*}
\sum_{n\ge j>l>1}\del_j(D_1+D_2)\cdot\del_l(D_1+D_2)=\sum_{m\ge j>l>1}\del_j(D_1)\cdot\del_l(D_1)\\
+\sum_{n\ge j>l\ge m}\del_j(D_2)\cdot\del_l(D_2)+\sum_{n\ge j\ge m}\sum_{m\ge l>1}\del_j(D_2)\cdot\del_l(D_1).
\end{align*}
After substituting all these, all that remains to be proved is that $$\del_1(D_2)\cdot\del D_1+\del D_2\cdot\del_1(D_1)+\sum_{n\ge j\ge m}\sum_{m\ge l>1}\del_j(D_2)\cdot\del_l(D_1)=0.$$
 However, the expression in question is simply $\del D_2\cdot\del D_1$ and hence is zero, thereby concluding the proof.
\end{proof}

Now, we try to understand how $\mu$ and $\iota$ behave under an isotopy of the $\eta$ curves. Let $\phi\in\pi_2(p^{1,2},\ldots,p^{n,1})$ be a Whitney $n$-gon. Let us do a small isotopy on $\eta^1$ curves to get $\wt{\eta^1}$ such that the $\eta^1$ curves stay disjoint throughout, $\wt{\eta^1}$ curves are transverse to the $\eta^i$ curves (however, the $\eta^1$ curves do not have to remain transverse to the $\eta^i$ curves throughout the isotopy) and the isotopy is constant in a neighborhood of the coordinates of $p^{1,2},\ldots,p^{n,1}$. The Whitney $n$-gon $\phi$ in the old Heegaard diagram $\mc{H}$ naturally gives rise to a Whitney $n$-gon $\wt{\phi}$ in the new Heegaard diagram $\wt{\mc{H}}$ joining the same $n$ points $p^{1,2},\ldots,p^{n,1}$. Since during the isotopy, the torus $T_1$ stayed disjoint from the fat diagonal, we get $\iota(\phi)=\iota(\wt{\phi})$; and since during the isotopy, the torus $T_1$ was untouched near $p^{1,2},\ldots,p^{n,1}$, using the alternate definition of the Maslov index, we get $\mu(\phi)=\mu(\wt{\phi})$. Once more a similar result holds for domains in the Heegaard diagram.

\begin{thm}\label{thm:isotopy}
Let $D\in\mc{D}(p^{1,2},\ldots,p^{n,1})$ be a domain. After an isotopy of the $\eta^1$ curves which is constant in a neighborhood of the coordinates of $p^{1,2},\ldots,p^{n,1}$, let $\wt{D}$ be the induced domain in the new Heegaard diagram joining the same $n$ points $p^{1,2},\ldots,p^{n,1}$. Then $\mu(D)=\mu(\wt{D})$ and $\iota(D)=\iota(\wt{D})$.
\end{thm}

\begin{proof}
The isotopy can be thought of as a finite sequence of steps, where at each step, exactly one of the two moves illustrated in Figure \ref{fig:moves} happens. The thin line denotes an $\eta^1$ curve, the thick line denotes an $\eta^u$ curve and the thick dotted line denotes an $\eta^v$ curve with $u,v\ne 1$ and $u\ne v$. The coefficients of the domains $D$ and $\wt{D}$ are shown. 

\begin{figure}[ht]
\psfrag{a1}{$a$}
\psfrag{a2}{$a$}
\psfrag{b1}{$b$}
\psfrag{b2}{$b$}
\psfrag{b3}{$b$}
\psfrag{c1}{$c$}
\psfrag{c2}{$c$}
\psfrag{d}{$a+c-b$}
\psfrag{u1}{$a+b+c$}
\psfrag{v1}{$a+b+c$}
\psfrag{u2}{$a+b$}
\psfrag{v2}{$a+b$}
\psfrag{u3}{$a$}
\psfrag{v3}{$d$}
\psfrag{u4}{$d-c$}
\psfrag{v4}{$d-c$}
\psfrag{u5}{$d-b-c$}
\psfrag{v5}{$d-b-c$}
\psfrag{u6}{$a+c$}
\psfrag{v6}{$a+c$}
\psfrag{u7}{$d-b$}
\psfrag{v7}{$d-b$}
\psfrag{x}{$\eta^1$}
\psfrag{y}{$\eta^u$}
\psfrag{z}{$\eta^v$}
\begin{center}
\includegraphics[width=280pt]{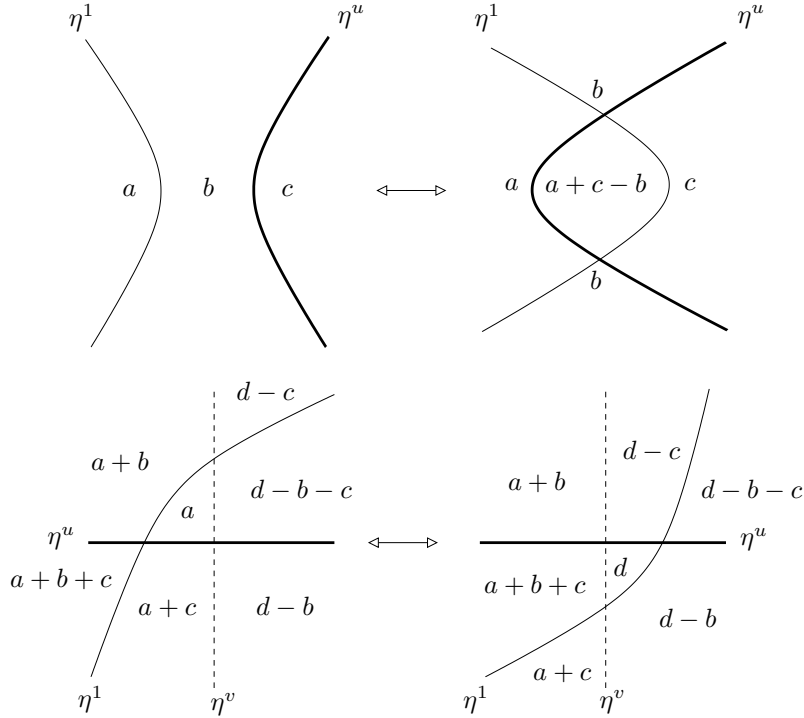}
\end{center}
\caption{Local coefficients of $D$ and $\wt{D}$}\label{fig:moves}
\end{figure}

Since $D$ and $\wt{D}$ have the same local coefficients near the
coordinates of $p^{1,2},\ldots,\allowbreak p^{n,1}$, the point measures $\mu_{p^{i,i+1}}$ do not change. From Figure \ref{fig:moves}, it is clear that $\del_i(D)\cdot\del_j(D)=\del_i(\wt{D})\cdot\del_j(\wt{D})$ for all values of $i,j$. Finally the identity $e(D)=e(\wt{D})$ can be verified in the following way.

Draw an extra circle around the local picture in each of the four cases, and choose a metric such that the circle is a geodesic which intersects all the $\eta$ curves in right angles. Since the Euler measure is additive, the Euler measure of the whole domain is the sum of the Euler measure of the domain lying outside the circle and the Euler measure of the domain lying inside the circle. Since the move is a local move, the Euler measure of the domain lying outside the circle does not change; and the identities $\frac{a}{2}+\frac{c}{2}=\frac{b}{4}+\frac{b}{4}+\frac{a+c-b}{2}$ and $\frac{a}{4}+\frac{d-c}{4}+\frac{d-b}{4}+\frac{a+b+c}{4}=\frac{d}{4}+
\frac{d-b-c}{4}+\frac{a+c}{4}+ \frac{a+b}{4}$ show that neither does the Euler measure of the domain lying inside the circle.
\end{proof}

We are almost done establishing all the relevant properties of $\mu$ and $\iota$. Let us end the section by doing a calculation for a special type of domain.

Assume that $n=3$, and let $D\in\mc{D}(p^{1,2},p^{2,3},p^{3,1})$. Furthermore assume that $D$ only has coefficients $0$ or $1$, and the union of the closure of the regions where $D$ has coefficient $1$ is a disjoint union of $k$ triangles. Then it is well known that there is Whitney $3$-gon $\phi\in\pi_2(p^{1,2},p^{2,3},p^{3,1})$ such that $D=D(\phi)$, and $\mu(\phi)=\iota(\phi)=0$. This is also true at the level of domains.

\begin{thm}\label{thm:testcase}
Let $D\in\mc{D}(p^{1,2},p^{2,3},p^{3,1})$ be a disjoint union of $k$ triangles. Then $\mu(D)=\iota(D)=0$.
\end{thm}

\begin{proof}
Number the $k$ triangles $\Del_1,\ldots,\Del_k$ arbitrarily, and assume without loss of generality that the vertices lying on the boundary of $\Del_i$ are $p^{1,2}_i$, $p^{2,3}_i$ and $p^{3,1}_i$. Therefore for all $i$, we have $e(\Del_i)=\frac{1}{4}$, $\del_3(\Del_i)\cdot\del_2(\Del_i)=-\frac{1}{4}$ and $\mu_{p^{3,1}_i}(\Del_i)=\mu_{p^{1,2}_i}(\Del_i)=\frac{1}{4}$. The Euler measure is additive, so $e(D)=\sum_i e(\Del_i)=\frac{k}{4}$; since the $k$ triangles are disjoint, we have $\del_3(D)\cdot\del_2(D)=\sum_i \del_3(\Del_i)\cdot\del_2(\Del_i)=-\frac{k}{4}$; and finally once more using the fact that the $k$ triangles are disjoint, we get $\mu_{p^{3,1}}(D)=\sum_i\mu_{p^{3,1}_i}(\Del_i)=\frac{k}{4}$ and $\mu_{p^{1,2}}(D)=\sum_i\mu_{p^{1,2}_i}(\Del_i)=\frac{k}{4}$. Therefore $\iota(D)=\mu_{p^{3,1}}(D)+\mu_{p^{1,2}}(D)+\del_3(D)\cdot\del_2(D)-e(D)=0$ and $\mu(D)=\iota(D)+2e(D)-\frac{k}{2}=0$.
\end{proof}

\section{$\mu(\phi)=\mu(D(\phi))$ and $\iota(\phi)=\iota(D(\phi))$}\label{sec:main}

We devote this section to proving the main theorem that the formulas for $\mu$ and $\iota$ for domains represent the the actual Maslov index and the intersection number in the symmetric product respectively.

\begin{thm}
For $n\geq 2$ and points $p^{i,i+1}\in P^{i,i+1}$, let $\phi\in\pi_2(p^{1,2},\ldots,p^{n,1})$ be a Whitney $n$-gon. Then $\mu(\phi)=\mu(D(\phi))$ and $\iota(\phi)=\iota(D(\phi))$.
\end{thm}

\begin{proof}
We prove this by an induction on $n$. The case for $n=2$ is a theorem of Robert Lipshitz \cite[Corollary 4.10]{RL}, which was also partially proved by Jacob Rasmussen in \cite[Theorem 9.1]{JR}. Therefore let us work with $n\ge 3$.

Without loss of generality, let us assume that the coordinates of $p^{1,2}$ and $p^{2,3}$ lying on $\eta^2_i$ are numbered $p^{1,2}_i$ and $p^{2,3}_i$ respectively. By renumbering the $\eta^1$ and the $\eta^3$ circles if necessary, let us also assume that $\eta^1_i$ passes through $p^{1,2}_i$ and $\eta^3_i$ passes through $p^{2,3}_i$.

Let $U_i$ be a small neighborhood of $\eta^2_i$. We have already fixed an orientation on $\Sigma$. Let us now fix some arbitrary orientations on the $\eta^2$ circles. We declare that each $\eta^2_i$ circle runs from west to east, therefore at each point in $U_i$, we have a notion of the directions east, west, north and south.

Let us choose a point $a_i\in U_i\cap\eta^3_i$ to the north of $p^{2,3}_i$, and let $\tau_i$ be an embedded path lying inside $U_i\setminus\eta^2_i$, joining $a_i$ to a point on $\eta^1_i$, such that $\tau_i$ either throughout travels northeastwards or throughout travels northwestwards, and $\tau_i$ is homotopic to $-\del_2(D(\phi))$ inside $U_i$ relative to the $\eta^1$ and the $\eta^3$ circles.

Now we will do a finger-move isotopy on the $\eta^3$ curves. Make a finger at $\eta^3_i$ near $a_i$, and push it along $\tau_i$ all the way until it intersects $\eta^1_i$ near the endpoint of $\tau_i$. We respect the following usual conventions of a finger move as specified in \cite{SSJW}. Whenever the finger encounters any $\eta^3$ curve, push it along with the finger, and whenever it encounters any other $\eta$ curve, make it intersect the finger. Let $q_i$ be the southern intersection of the innermost finger with $\eta^1_i$ near the endpoint of $\tau_i$. Choose a point $b_i$ just to the east of $q_i$ on the new $\eta$ curve, make a finger there, and push it all the way south until it hits $\eta^2_i$. Let $r_i$ be the western intersection of the new innermost finger with $\eta^2_i$. Let $\eta^{3'}_i$ be the new $\eta^3_i$ curve, and let $\eta^{3'}=\{\eta^{3'}_1,\ldots,\eta^{3'}_k\}$. We call $\mc{H}=(\Sigma,\eta^1,\eta^2,\eta^3,\ldots,\eta^n)$ the old Heegaard diagram and $\mc{H}'=(\Sigma,\eta^1,\eta^2,\eta^{3'},\ldots,\eta^n)$ the new Heegaard diagram. Let $T_{3'}\subset Sym^k(\Sigma)$ be the torus corresponding to $\eta^{3'}$ and let $P^{i,3'}=T_i\cap T_{3'}$. The isotopy is illustrated in Figure \ref{fig:main}; the thin dotted lines represent $\tau_i$, the thin solid lines represent $\eta^1_i$, the thick dotted lines represent $\eta^2_i$ and the thick solid lines represent either $\eta^3_i$ or $\eta^{3'}_i$ or the intermediate curve. The points $p^{1,2}_i$, $p^{2,3}_i$, $a_i$, $b_i$, $q_i$ and $r_i$ are shown.

\begin{figure}
\psfrag{a}{$\eta^1_i$}
\psfrag{b}{$\eta^2_i$}
\psfrag{c}{$\eta^3_i$}
\psfrag{cc}{$\eta^{3'}_i$}
\psfrag{aa}{$a_i$}
\psfrag{bb}{$b_i$}
\psfrag{r}{$r_i$}
\psfrag{q}{$q_i$}
\psfrag{p23}{$p^{2,3}_i$}
\psfrag{p23x}{$p^{2,3}_i=x^i_0$}
\psfrag{p12}{$p^{1,2}_i$}
\psfrag{x}{$x^i_1$}
\psfrag{xx}{$x^i_2$}
\psfrag{t}{$\tau_i$}
\begin{center}
\includegraphics[width=370pt]{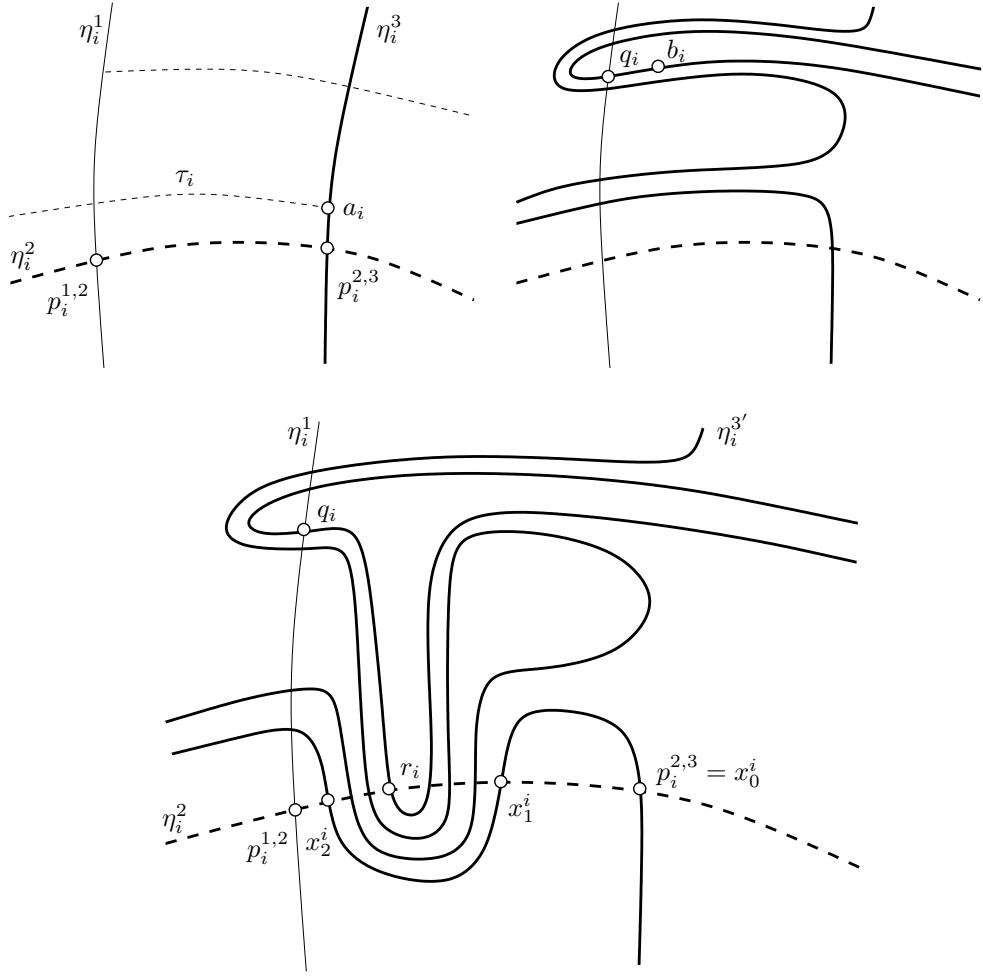}
\end{center}
\caption{The isotopy of the $\eta^3$ circles in $U_i$}\label{fig:main}
\end{figure}

Note that the isotopy of the $\eta^3$ circles is constant in a neighborhood of the coordinates of $p^{i,i+1}$. Therefore the Whitney $n$-gon $\phi$ in $\mathcal{H}$ gives rise to a Whitney $n$-gon $\phi '$ in $\mathcal{H}'$ also joining the points $p^{1,2},\ldots,p^{n,1}$, such that $\mu(\phi)=\mu(\phi')$, $\iota(\phi)=\iota(\phi')$, $\mu(D(\phi))=\mu(D(\phi'))$ and $\iota(D(\phi))=\iota(D(\phi'))$.

Let us travel on the $\eta^{3'}_i$ curve from $p^{2,3}_i$ to $r_i$, and let $x^i_0,x^i_1,\ldots,x^i_{n_i}$ be in order, the points of intersection with $\eta^2_i$ that we encounter on the way, such that $x^i_0=p^{2,3}_i$ and $x^i_{n_i}=r_i$. Let $q=(q_1,\ldots,q_k)\in P^{1,3'}$, $r=(r_1,\ldots,r_k)\in P^{2,3'}$ and $y^i_j=(r_1,\dots,r_{i-1},x^i_j,p^{2,3}_{i+1},\dots,p^{2,3}_k)\in P^{2,3'}$. Note that $y^1_0=p^{2,3}$, $y^{i+1}_0=y^i_{n_i}$ and $y^k_{n_k}=r$.

For each $i,j$, there is an obvious Whitney $2$-gon $u^i_j\in\pi_2(y^i_j,y^i_{j+1})$ such that $D(u^i_j)$ either only has coefficients $0$ or $1$, or only has coefficients $0$ and $-1$, and the closure of the union of the regions where $D(u^i_j)$ has non-zero coefficients is a bigon supported inside $U_i$. Either by a direct computation, or from Lipshitz' formula \cite[Corollary 4.10]{RL}, we know that $\mu(u^i_j)=\mu(D(u^i_j))$ and $\iota(u^i_j)=\iota(D(u^i_j))$ for all $i,j$. Therefore the multiplication operation for Whitney $2$-gons, produces a Whitney $2$-gon $u=*_{i=1}^k *_{j=0}^{n_i-1}u^i_j\in\pi_2(p^{2,3},r)$ such that $\mu(u)=\mu(D(u))$ and $\iota(u)=\iota(D(u))$.

There is also an obvious Whitney $3$-gon $v\in\pi_2(r,q,p^{1,2})$ such that $D(v)$ only has coefficients $0$ and $1$, and the closure of the union of the regions where $D(v)$ has non-zero coefficients is a disjoint union of $k$ triangles, one in each $U_i$. Therefore by Theorem \ref{thm:testcase}, $\mu(v)=\mu(D(v))$ and $\iota(v)=\iota(D(v))$. Theorem \ref{thm:addition} now implies that $\mu(u*v)=\mu(D(u*v))$ and $\iota(u*v)=\iota(D(u*v))$.

By construction, we have $\del_2(D(u*v))=\del_2(D(\phi))=\del_2(D(\phi'))$. The Whitney $3$-gon $(u*v)$ and the Whitney $n$-gon $\phi'$ can be represented by two maps from the unit disk $D$ to $Sym^k(\Sigma)$ such that there is a fixed arc $s_2\subset\del D$ which maps to $T_2$ and whose boundary points map to $p^{1,2}$ and $p^{2,3}$. Since $\del_2(D(u*v))=\del_2(D(\phi'))$, the two images of $s_2$ are homotopic in $T_2$ relative the endpoints. Therefore there exists a Whitney $(n-1)$-gon $\psi\in\pi_2(q,\ldots,p^{n,1})$ such that $(u*v)*\psi=\phi'$. Theorem \ref{thm:addition} then implies that $\mu(\phi')=\mu(\psi)+\mu(u*v)$, $\mu(D(\phi'))=\mu(D(\psi))+\mu(D(u*v))$, $\iota(\phi')=\iota(\psi)+\iota(u*v)$ and $\iota(D(\phi'))=\iota(D(\psi))+\iota(D(u*v))$.

We now fit all the pieces to complete the proof. The induction hypothesis gives us the starting block $\mu(\psi)=\mu(D(\psi))$ and $\iota(\psi)=\iota(D(\psi))$. The rest follows formally; $\mu(\phi)=\mu(\phi')=\mu(\psi)+\mu(u*v)=\mu(D(\psi))+\mu(D(u*v))=\mu(D(\phi'))=\mu(D(\phi))$; similarly $\iota(\phi)=\iota(\phi')=\iota(\psi)+\iota(u*v)=\iota(D(\psi))+\iota(D(u*v))=\iota(D(\phi'))=\iota(D(\phi))$.
\end{proof}

\section{Some applications}\label{sec:applications}
In this section we present some well known examples and applications. These require some concepts in addition to the ones that we discussed in Section \ref{sec:setting}. To avoid unnecessary cluttering in that section, we introduce these new concepts in this section as and when we need them. 

\subsection{Absolute grading formula}

Let $\mc{H}=(\Sigma,\eta^1,\eta^2,\eta^3, w)$ be a Heegaard diagram with $n=3$ and $k=1$, and let $\phi,\phi'\in\pi_2(p^{1,2},p^{2,3},p^{3,1})$ be two Whitney $3$-gons joining $p^{1,2}$, $p^{2,3}$ and $p^{3,1}$. Let the $2$-chain $D(\phi')-D(\phi)$ be denoted by $P$. The absolute grading formula in \cite[Formula (12)]{POZSz4manifolds} suggests that $\mu(\phi')-\mu(\phi)=2n_w(P)+\frac{c_1(\mathfrak{s}_w(\phi'))^2-
c_1(\mathfrak{s}_w(\phi))^2}{4}$. Here $n_w(P)$ is the coefficient of the $2$-chain $P$ at the region containing the basepoint $w$; given a $Spin^C$ structure $\mf{s}$ on the $4$-manifold $W_{\eta^1,\eta^2,\eta^3}$, $c_1(\mf{s})$ denotes its first Chern class; and $\mf{s}_w(\phi)$ is the $Spin^C$ structure on $W_{\eta^1,\eta^2,\eta^3}$ coming from a Whitney $3$-gon $\phi$ with respect to the basepoint $w$, see \cite[Section 2.2]{POZSz4manifolds}. In this subsection, we  give a direct verification of the above formula using the formula for the Maslov index $\mu$.

Following \cite[Section 2.2]{POZSz4manifolds}, let $H(P)\in H_2(W_{\eta^1,\eta^2,\eta^3})$ be the homology class corresponding to $P$, and let $PD(H(P))$ denote its \Poincare{} dual.
The right hand side then simplifies to $2n_w(P)+\langle c_1(\mathfrak{s}_w(\phi)), H(P)\rangle +PD(H(P))^2$.
However, yet another formula from \cite[Proposition 6.3]{POZSz4manifolds} implies that $\langle c_1(\mathfrak{s}_w(\phi)), P\rangle=e(P)+\#(\partial P)-2n_w(P)+2\sigma(\phi,P)$. Here $\#(\del P)$ denotes the number of components in $\del P$ counted with multiplicity; and $\sigma(\phi,P)$ denotes the dual spider number defined as follows. The Whitney $3$-gon is represented by a map from the unit disk $D^2$ to the symmetric product. There are three arcs $s_1$, $s_2$ and $s_3$ on $\del D^2$ which map to the tori $T_1$, $T_2$ and $T_3$ respectively. Choose a point $x$ in the interior of $D^2$ and three arcs $e_1$, $e_2$ and $e_3$, such that $e_i$ joins $x$ to a point in the interior of $s_i$. Assume that the image of each $e_i$ is disjoint from the fat diagonal. Therefore, if the image of $x$ is $(y_1,\ldots,y_k)$, then the image of $e_i$ can be represented by $k$ arcs on $\Sigma$, $f^i_1,\ldots,f^i_k$, such that $f^i_j$ joins $y_j$ to a point on the $\eta^i$ circles. Let $\del'_i(P)$ be a small outward translate of the circles in $\del_i(P)$. Then the dual spider number is defined as $\sigma(\phi,P)=\sum_j(n_{y_j}(P)+\sum_if^i_j\cdot\del'_i(P))$. Therefore to complete our verification, we only need to check that $$\mu(D(\phi)+P)-\mu(D(\phi))=e(P)+\#(\del P)+2\sigma(\phi,P)+PD(H(P))^2.$$

We have already seen that our formula for $\mu$ is cyclically
symmetric in $p^{1,2},p^{2,3},\allowbreak p^{3,1}$, so we use a more symmetric version where $$\mu(D)=e(D)+\frac{2}{3}\sum_i\mu_{p^{i,i+1}}(D)+\frac{1}{3}\sum_i\del_{i+1}(D)\cdot\del_i (D).$$
 Therefore, the left hand side simplifies to  
\begin{align*}
\frac{1}{3}\sum_i (\del_{i+1}(D(\phi)).\del_i(P) +\del_{i+1}(P).\del_i(D(\phi))+\del_{i+1}(P).\del_i(P))\\
+e(P)+\frac{2}{3}\sum_i\mu_{p^{i,i+1}}(P).
\end{align*} 

It follows from the definition of the cohomology class $PD(H(P))$ that $PD(H(P))^2\allowbreak=\del_{i+1}(P)\cdot\del_i(P)$ for all $i$. Therefore, we only need to show that
\begin{align*}
\frac{2}{3}\sum_i\mu_{p^{i,i+1}}(P)+\frac{1}{3}\sum_i (\del_{i+1}(D(\phi)).\del_i(P) +\del_{i+1}(P).\del_i(D(\phi)))\\=\#(\partial P)+2\sigma(\phi,P).
\end{align*}

In order to calculate $\sigma(\phi,P)$, choose the interior point $x$ close to the point $t_2=s_2\cap s_3$, have the arcs $e_2$ and $e_3$  be small and supported near $t_2$, and ensure that the arc $e_1$ is supported near $s_3$ and runs parallel to it. Then $$\sum_j(n_{y_j}(P)+f^2_j\cdot\del'_2(P)+f^3_j\cdot\del'_3(P))=\mu_{p^{2,3}}(P)-\frac{1}{2}(\#(\del_2P)+\#(\del_3P)),$$ $$\text{and } \sum_j(f^1_j\cdot\del'_1(P))=\del_1(P)\cdot\del_3(D(\phi))-\frac{1}{2}(\#(\del_1P)).$$ Therefore we get $$\sigma(\phi,P)=\mu_{p^{2,3}}(P)-\frac{1}{2}(\#(\del P))+\del_1(P)\cdot\del_3(D(\phi)).$$ 

Instead we could have chosen the arc $e_1$ to be parallel to $s_2$, and we would have got $$\sigma(\phi,P)=\mu_{p^{2,3}}(P)-\frac{1}{2}(\#(\del P))+\del_2(D(\phi))\cdot\del_1(P).$$ 
Adding, we get, 
$$2\sigma(\phi,P)=2\mu_{p^{2,3}}(P)-\#(\del P)+\del_1(P)\cdot\del_3(D(\phi))+\del_2(D(\phi))\cdot\del_1(P).$$ However, by choosing the point $x$ near $t_1$ or $t_3$, we get two similar expressions for the dual spider number. After taking averages, we obtain our required identity
\begin{align*}
2\sigma(\phi,P)=\frac{1}{3}\sum_i (\del_{i+1}(D(\phi)).\del_i(P) +\del_{i+1}(P).\del_i(D(\phi)))\\+\frac{2}{3}\sum_i\mu_{p^{i,i+1}}(P)-\#(\partial P).
\end{align*}

While we are on the topic of absolute gradings, it might be interesting to note that a combinatorial formula for Maslov index of triangles allows us to compute the absolute Maslov grading on $\wh{HF}(Y,\mf{t})$, the hat version of the Heegaard Floer homology of a three-manifold $Y$ in a torsion $Spin^C$ structure $\mf{t}$. This follows from the proof of Theorem 7.1 in \cite{POZSz4manifolds}. 
There is a Heegaard diagram $\mc{H}=(\Sigma, \eta^1,\eta^2,\eta^3,w)$ with $n=3$ and $k=1$, such that $Y_{\eta^1,\eta^2}=S^3$, $Y_{\eta^2,\eta^3}=\#^l(S^2\times S^1)$ for some $l$, and $Y_{\eta^1,\eta^3}=Y$. A Whitney $3$-gon $\psi\in\pi_2(p^{1,2},p^{2,3},p^{3,1})$ is chosen such that the $Spin^C$ structure associated to $\psi$ restricts to $\mf{t}$ on $Y$, and to the torsion $Spin^C$ structure $\mf{t}_0$ on $\#^l(S^2\times S^1)$. Then \cite[Formula (12)]{POZSz4manifolds} provides a combinatorial formula relating $\mu(\psi)$ and the absolute gradings of $p^{1,2},p^{2,3}$ and $p^{3,1}$. However, the Heegaard diagrams $(\Sigma,\eta^1,\eta^2)$ and $(\Sigma,\eta^2,\eta^3)$ can be modified by sequences of isotopies and handleslides until they become `standard' Heegaard diagrams for $S^3$ and $\#^l(S^2\times S^1)$ respectively. Thus the absolute gradings of $p^{1,2}$ and $p^{2,3}$ can be determined. Therefore, a combinatorial formula of $\mu(\psi)$ leads to an algorithm to compute the absolute Maslov grading of $p^{3,1}$. However, the Heegaard diagram $(\Sigma,\eta^1,\eta^3)$ can be converted by a sequence of isotopies and handleslides to a Heegaard diagram where $\wh{HF}(Y,\mf{t})$ can be computed combinatorially, see \cite{SSJW}. The absolute grading of $p^{3,1}$ allows us to assign absolute gradings to all the generators in the new Heegaard diagram that lie in the $Spin^C$ structure $\mf{t}$, and thereby completes the computation of the absolute Maslov grading on $\wh{HF}(Y,\mf{t})$.

\subsection{Domains supported on regions with non-negative Euler measure}

In a Heegaard diagram, a positive domain is a domain such that none of its coefficients are negative and at least one of its coefficients is positive. In the Heegaard Floer world, we are often interested in positive domains. This is because, if the complex structure on the symmetric product $Sym^k(\Sigma)$ is sufficiently close to the complex structure induced from one on $\Sigma$, and if the moduli space $\mc{M}(\phi)$ is non-empty, then the domain $D(\phi)$ is either the trivial domain or a positive domain.

We are also interested in domains whose only non-zero coefficients lie on regions that are topologically disks with non-negative Euler measure, or in other words $m$-sided regions with $m\leq 4$. Such domains come up naturally when we work with Heegaard diagrams where any region that does not contain any basepoint $w_j$ is either a $2$-sided or a $3$-sided or a $4$-sided region, and we are only interested in domains that avoid the basepoints. Such Heegaard diagrams are usually called nice Heegaard diagrams, but nice Heegaard diagrams are usually so complicated that the terminology is at best a misnomer.

Let $\mc{H}=(\Sigma,\eta^1,\ldots,\eta^n,w)$ be a nice Heegaard diagram. Nice Heegaard diagrams were first studied in \cite[Definition 3.1]{SSJW}, where for $n=2$, it was proved that the computation of the moduli space $\mc{M}(\phi)$ for a Whitney $2$-gon $\phi$ is combinatorial in a nice Heegaard diagram if $\mu(\phi)=1$ and $D(\phi)$ avoids the basepoints in $w$. Robert Lipshitz, Ciprian Manolescu and Jiajun Wang studied nice Heegaard diagrams for $n=3$ in \cite{RLCMJW}, where they independently proved a theorem that we are going to prove shortly, that the computation of the moduli space $\mc{M}(\phi)$ for a Whitney $3$-gon $\phi$ is also combinatorial in a nice Heegaard diagram if $\mu(\phi)=0$ and $D(\phi)$ avoids the basepoints in $w$.

\begin{thm}\label{thm:nice}
Let $\mc{H}=(\Sigma,\eta^1,\eta^2,\eta^3)$ be a Heegaard diagram, and
let $D\in\mc{D}(p^{1,2},\allowbreak p^{2,3},\allowbreak p^{3,1})$ be a positive domain such that $D$ is supported on $2$-sided, $3$-sided and $4$-sided regions, and $\iota(D)\geq 0$. Then $\mu(D)\geq 0$, and equality holds only if $\iota(D)=0$.
\end{thm}

\begin{proof}
First, observe that $\mu(D)=\iota(D)+2e(D)-\frac{k}{2}\geq 2e(D)-\frac{k}{2}$. Next recall that $D$ is a positive domain supported on regions with non-negative Euler measure, therefore $e(D)=\frac{1}{4}($number of $3$-sided regions in the $2$-chain $D)+\frac{1}{2}($number of $2$-sided regions in the $2$-chain $D)$. 

Given a $2$-chain $D'$, regard $\del(\del_1 D')_{|\eta^2}$ as a formal sum of points, and let $s(D')$ be the sum of the coefficients of the points in that formal sum. Clearly $s$ is an additive function and $s(D)=k$. However if $R$ is a $2$-sided or a $4$-sided region then $s(R)=0$, and if $R$ is a $3$-sided region, then $s(R)=\pm 1$. Therefore, counted with multiplicities, $D$ must contain at least $k$ $3$-sided regions, which in turn proves that $e(D)\geq \frac{k}{4}$ and $\mu(D)\geq 0$.

If $\mu(D)=0$, then equality must hold throughout. Therefore $\iota(D)=0$ and $e(D)=\frac{k}{4}$. This then implies that $D$ contains some number of $4$-sided regions, exactly $k$ $3$-sided regions, each with $s=1$, and no $2$-sided region.
\end{proof}

The following theorem shows that the moduli space $\mc{M}(\phi)$ can be determined for Whitney $3$-gons with Maslov index $0$. We will use Lipshitz' reformulation \cite{RL} of Heegaard Floer theory, so let us briefly mention it now.

Recall that a Whitney $3$-gon $\phi\in\pi_2(p^{1,2},p^{2,3},p^{3,1})$ can be represented by a map from the unit disk $D^2$ to the symmetric product $Sym^k(\Sigma)$. There are three points $t_1$, $t_2$ and $t_3$ on $\del D^2$ and three positively oriented arcs $s_1$, $s_2$ and $s_3$ also on $\del D^2$, such that $s_i$ is disjoint from $t_{i+1}$ and it joins $t_{i-1}$ to $t_i$. The map from $D^2$ to $Sym^k(\Sigma)$ is required to map $t_i$ to $p^{i,i+1}$ and $s_i$ to $T_i$ for all $i$. If such a map can be represented by a holomorphic map (with respect to the complex structure on $Sym^k(\Sigma)$ induced from one on $\Sigma$), i.e. if $\mc{M}(\phi)\neq\varnothing$, then such a map can also be represented by a Lipshitz map, as described below.

Given a domain $D\in\mc{D}(p^{1,2},p^{2,3},p^{3,1})$, a Lipshitz map is a pair $(F,u)$, where $F$ is a surface with a complex structure and $u$ is a holomorphic embedding of $F$ into $D^2\times\Sigma$ such that the following conditions are satisfied; if $p_1$ is the first projection map, then $p_1u$ is a $k$-sheeted branched covering, with all the branch points lying in the interior of $D^2$; if $p_2$ is the second projection map, then the degree of the map $p_2u$ at a region is the coefficient of $D$ at that region; and finally, we have for all $i$, $p_2u(p_1u)^{-1}(t_i)=p^{i,i+1}$ and $p_2u(p_1u)^{-1}(s_i)\subseteq\eta^i$. 

It turns out that if a domain $D$ can be represented by a Lipshitz map, then there is a Whitney $3$-gon $\phi$ such that $D=D(\phi)$, and in that case, the number of branch points of $p_1u$ equals $\iota(\phi)$. Another important lemma is that, given a Whitney $3$-gon $\phi$, there is a bijection between $\mc{M}(\phi)$ and the set of all Lipshitz maps representing the domain $D(\phi)$.  With these facts in mind, we are all set to proceed to the next theorem.

\begin{thm}
Let $\mc{H}=(\Sigma,\eta^1,\eta^2,\eta^3,w)$ be a nice Heegaard diagram, and let $D\in\mc{D}(p^{1,2},p^{2,3},p^{3,1})$ be a domain which avoids the basepoints such that $\mu(D)=0$. We choose the complex structure on $Sym^k(\Sigma)$ induced from one on $\Sigma$. 
Then the following three statements are equivalent. 

(1) The domain $D$ is equal to $D(\phi)$ for some Whitney $3$-gon $\phi$ and $\mc{M}(\phi)$ has one point.

(2) The domain $D$ is equal to $D(\phi)$ for some Whitney $3$-gon $\phi$ and $\mc{M}(\phi)$ is non-empty.

(3) The intersection number $\iota(D)=0$ and the $2$-chain $D$ can be represented by a (not necessarily disjoint) union of $k$ embedded triangles, such that the $3k$ sides of the $k$ triangles all lie on different $\eta$ circles.
\end{thm}

\begin{proof}
The equivalence of the three statements follows from the following three implications.

\textbf{(1)$\Rightarrow$(2)} This is quite obvious.

\textbf{(2)$\Rightarrow$(3)} Since the moduli space $\mc{M}(\phi)$ is non-empty, the intersection number $\iota(D)$ is non-negative. Furthermore there is some Lipshitz map $(F,u)$ representing $D$. The image of $F$ under the map $p_2u$ at the $2$-chain level is $D$, therefore $D$ is a positive domain. Theorem \ref{thm:nice} applies and we get that $\iota(D)=0$ and $D$ contains some number of $4$-sided regions, exactly $k$ $3$-sided regions and no  $2$-sided region. The surface $F$ is a $k$-sheeted branched cover over the disk with $\iota(D)=0$ branch points, therefore $F$ is a disjoint union of $k$ disks $\Delta_1,\ldots,\Delta_k$, each with three marked points on its boundary.

Let $u_i=u_{|\Delta_i}$, let $p^{1,2}_i$, $p^{2,3}_i$ and $p^{3,1}_i$ be the images of the three marked points on $\del\Delta_i$, and let $\eta^1_i$, $\eta^2_i$ and $\eta^3_i$ be the $\eta$ circles passing through those three points. Therefore $\del \Delta_i$ maps to $\cup_j\eta^j_i$; hence $u_i$ is a Whitney $3$-gon in the Heegaard diagram $\mc{H}_i=(\Sigma,\eta^1_i,\eta^2_i,\eta^3_i)$ joining $p^{1,2}_i$, $p^{2,3}_i$ and $p^{3,1}_i$. The image of $u_i$ is a positive domain $D_i$ whose support lies on $3$-sided and $4$-sided regions and $\iota(u_i)=0$. In fact, since $s(u_i)=1$, the $2$-chain  $D_i$ contains at least one $3$-sided region. However the total number of $3$-sided regions (counted with multiplicities) in $D=\sum_i D_i$ is $k$, therefore $D_i$ contains exactly one $3$-sided region and some number of $4$-sided regions. This in turn implies that the Euler measure of $D_i$ is $\frac{1}{4}$, therefore the Maslov index $\mu(u_i)$ is $\iota(u_i)+2e(D_i)-\frac{1}{2}=0$. 

\begin{figure}
\psfrag{a}{$\eta^1$}
\psfrag{b}{$\eta^2$}
\psfrag{c}{$\eta^3$}
\begin{center}
\includegraphics[width=170pt]{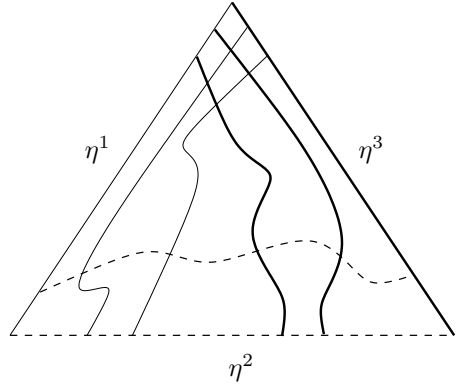}
\end{center}
\caption{Induced tiling on $\Delta_i$}\label{fig:tiling}
\end{figure}

We would like to show that each $D_i$ is an embedded triangle, or in other words, $D_i$ only has coefficients $0$ and $1$, and the closure of the union of the regions where $D_i$ has non-zero coefficients is a triangle. This would complete the second part of the theorem. Therefore we only need to prove that the map $u_i$ is a diffeomorphism.

Observe that the Euler measure of $D_i$ is $\frac{1}{4}$ and the Euler measure of $\Delta_i$, a disk with three marked points on its boundary, is also $\frac{1}{4}$. Therefore the map $u_i$ is an unbranched map, or in other words, a local diffeomorphism. Consider the preimages of the $\eta$ curves in $\Delta_i$. The preimage of each curve is an $1$-manifold, and the induced tiling on $\Delta_i$ has one $3$-sided region and some number of $4$-sided regions. It is not very hard then, to see that the induced tiling must look somewhat like Figure \ref{fig:tiling}. The preimages of the $\eta^1$ curves are denoted by thin lines; the preimages of the $\eta^2$ curves are denoted by thick dotted lines; and the preimages of the $\eta^3$ curves are denoted by thick solid lines. By an abuse of notation, the preimage of an $\eta^i$ curve is also called an $\eta^i$ curve. The intersections among these $\eta^i$ curves in $\Delta_i$ are called vertices.

Assume if possible, that $u_i$ is not injective. Therefore there are two distinct vertices $p$ and $q$ in $\Delta_i$, such that $u_i(p)=u_i(q)$. The vertices of $\Delta_i$ are naturally grouped into three groups based on whether they lie on $\eta^1$ and $\eta^2$ curves, $\eta^2$ and $\eta^3$ curves, or $\eta^3$ and $\eta^1$ curves. It is clear that $p$ and $q$ have to belong to the same group, and let us assume without loss of generality that they belong to the group which lies on the $\eta^1$ and $\eta^2$ curves.

Since $u_i$ is a local diffeomorphism, there exists a direction on the $\eta^1$  curve through $p$ and there exists a direction on the preimage of $\eta^1$ curve through $q$, such that $p$ and $q$ can be moved along these directions on the $\eta^1$ curves while ensuring $u_i(p)=u_i(q)$. The points $p$ and $q$ move along parallel (if not the same) curves, and therefore we can talk about whether they are moving in the same direction or in opposite directions. Furthermore notice that $p$ encounters a vertex on its way when and only when $q$ encounters a vertex on its way. Finally observe that since $u_i$ is an unbranched map, the points $p$ and $q$ remain disjoint. A similar statement holds for the $\eta^2$ curves. Since the map $u_i$ is orientation preserving, the points $p$ and $q$ move in the same direction along the $\eta^1$ curves if and only if they move in the same direction along the $\eta^2$ curves. There are two natural cases.

\emph{Case 1: The points $p$ and $q$ move in the same direction.}

Let us move $p$ along the $\eta^1$ curves towards the $\eta^3$ curves. Therefore $q$ also moves along the $\eta^1$ curves towards the $\eta^3$ curves. Since the condition $u_i(p)=u_i(q)$ holds true, $p$ reaches the first $\eta^3$ curve exactly when $q$ reaches the first $\eta^3$ curve. Therefore $p$ and $q$ must have crossed the same number of $\eta^2$ curves along the way, and hence $p$ and $q$ lie on the same $\eta^2$ curve. A similar argument shows that they lie on the same $\eta^1$ curve, thereby proving $p=q$.

\emph{Case 2: The points $p$ and $q$ move in opposite directions.}

Consider the rectangle $S$ in $\Delta_i$ that has $p$ and $q$ as its diametrically opposite corners. It is possible that the rectangle is degenerate, but that presents no problem. Let us move $p$ towards the center of $S$. Since $q$ moves in the opposite direction, $q$ also moves towards the center of $S$. If $S$ is tiled by an even number of $4$-sided regions (or as a special case, if $S$ is degenerate), then the center of $S$ lies on its $1$-skeleton. Therefore eventually $p$ and $q$ will hit the same point, which is a contradiction. On the other hand, if $S$ is tiled by an odd number of $4$-sided regions, then the center of $S$ lies in the interior of a $4$-sided region $S'$. Then $p$ and $q$ can be  moved such that they lie on two diametrically opposite corners of $S'$. Therefore $u_{i|S'}$ has degree at least $2$, which is a contradiction to the assumption that $u_i$ is a local diffeomorphism.

\textbf{(3)$\Rightarrow$(1)} We show that there exist one and only one Lipshitz map $(F,u)$ representing the domain $D$. Due to Lipshitz' reformulation, this is enough to establish that $\mc{M}(\phi)$ consists of exactly one point. 

Since $\iota(D)=0$, $F$ must be a disjoint union of $k$ disks $\Delta_1,\ldots,\Delta_k$, each with $3$ marked points on its boundary. Each such disk $\Delta_i$ admits a unique holomorphic map to the unit disk $D^2$, which maps the three marked points on $\del\Delta_i$ to the three marked points on $\del D^2$. The domain $D$ is a union of $k$ embedded triangles, such that the $3k$ sides on the boundary of the $k$ triangles all lie on different $\eta$ circles. Therefore $D$ can be realized as the image of $F=\cup_i\Delta_i$ under some map, and for each value of $i$ and for any of the $k$ triangles, there is a unique holomorphic map from $\Delta_i$ to that triangle, which sends the marked points on $\del\Delta_i$ to the vertices of the triangle.

Therefore we see that given such a domain $D$, there can be at most one Lipshitz map $(F,u)$ representing $D$. In fact, we have almost constructed the unique Lipshitz map. There is a surface $F$ which admits a holomorphic $k$-sheeted branched cover over the unit disk $D^2$ with $\iota(D)$ branch points, and also a holomorphic map to $\Sigma$ such that the image at the $2$-chain level is the domain $D$. We only need to show that the induced map $u:F\rightarrow D^2\times\Sigma$ is an embedding. Assume that the surface $F$ lying inside $D^2\times\Sigma$ has $d$ double points. We can modify the surface $F$ near the $d$ double points to obtain a new surface $F'$ which is embedded in $D^2\times\Sigma$, and $\chi(F')=\chi(F)-2d$. Therefore the new pair $(F',u')$ is a Lipshitz map representing the domain $D$, and hence the map $p_1u'$ must have $\iota(D)$ branch points. However the number of branch points of $p_1u'$ is $\iota(D)+2d$, therefore $d=0$. 
\end{proof}

\bibliographystyle{amsalpha}

\bibliography{maslov}

\end{document}